\documentclass[12pt]{article}

\def\ee{\epsilon}
\def\aa{\alpha}

\usepackage{latexsym}

\def\be{\begin{equation}}
\def\ee{\end{equation}}

\newcommand{\hatl}{\hat L}
\newcommand{\D}{\hat D}
\newcommand{\X}{\hat X}

\newtheorem{theorem}{Theorem}

\usepackage{graphicx}

\begin{document}

\date{}

\title{Exact solutions of hyperbolic systems of kinetic equations.
Application to Verhulst model with random perturbation\thanks{This
paper was written with partial financial support from the RFBR
grant 06-01-00814 and the DFG Research Unit 565 `Polyhedral
Surfaces'' (TU-Berlin).}}

\author{E.I.~Ganzha, V.M.~Loginov,
S.P.~Tsarev}

   \maketitle
\begin{center}
Krasnoyarsk State Pedagogical University\\
 89 Lebedevoi, 660049
Krasnoyarsk, Russia\\[1ex]
and \\[1ex]
 Department of Mathematics \\
 Technische Universit\"at Berlin \\
Berlin, Germany \\[1ex] 
e-mails: \\
\texttt{loginov@imfi.kspu.ru}\\
\texttt{tsarev@math.tu-berlin.de} \ 
\texttt{tsarev@newmail.ru}\\
\end{center}

\medskip

\begin{abstract}

 For hyperbolic first-order systems of linear partial differential
equations (master equations), appearing in description of kinetic
processes in physics, biology and chemistry we propose a new
procedure to obtain their complete closed-form non-stationary
solutions. The methods used include the classical Laplace cascade
method as well as its recent generalizations for systems with more
than 2 equations and more than 2 independent variables. As an
example we present the complete non-stationary solution
(probability distribution) for Verhulst model driven by Markovian
coloured dichotomous noise.


\end{abstract}

\section{Introduction}

This paper is devoted to a novel application of methods of
explicit integration of hyperbolic linear systems of PDEs recently
developed in \cite{ts98,ts05,ts06} to an important class of
dynamical nonlinear systems driven by a coloured noise.

Modelling dynamical systems, in which stochastic behaviour is
involved, as a rule leads to nonlinear stochastic differential
equations for the dynamical variable or sets of dynamical
variables. The comprehensive statistical treatment of these
variables may be formulated with the aid of Liouville stochastic
equation for the probability distribution (see \cite{HL84}).

{\it Example 1.} In this paper we consider as the simplest example
the following one-dimensional dynamical system
\begin{equation}\label{Ver}
\dot{x} = p(x)+\alpha(t) q(x),
\end{equation}
where $x(t)$ is the dynamical variable, $p(x)$, $q(x)$ are given
functions of $x$, $\aa(t)$ is the random function with known
statistical characteristics. The model (\ref{Ver}) arises in
different applications (see for example \cite{HL84,Kampen} and
bibliography therein). An important application of this model
consists in study of noise-induced transitions in physics,
chemistry and biology. The functions $p(x)$, $q(x)$ are often
taken polynomial. For example, if we set $p(x)=p_1x+p_2x^2$,
$q(x)=q_2x^2$, $p_1>0$, $p_2<0$, $|p_2|>q_2>0$, then the equation
(\ref{Ver}) describes the population dynamics when resources
(nutrition) fluctuate (Verhulst model).
 In the following we will assume $\aa(t)$ to be
binary (dichotomic) noise $\alpha(t)=\pm 1$ with switching
frequency $2\nu>0$. As one can show (see \cite{ShL78,LL}), the
averages $W(x,t)=\langle \widetilde W(x,t)\rangle$ and
$W_1(x,t)=\langle \alpha(t)\widetilde W(x,t)\rangle$ for the
probability density $\widetilde W(x,t)$ in the space of possible
trajectories $x(t)$ of the ODE satisfy a system of the form
(\ref{7}) (also called ``master equations"):
\begin{equation}\label{WW}
   \left\lbrace  \begin{array}{l}
 W_t + \left(p(x)W\right)_x  + \left(q(x)W_1\right)_x =0, \\[0.5em]
 (W_1)_t + 2\nu W_1 +\left(p(x)W_1\right)_x  + \left(q(x)W\right)_x
 =0.
 \end{array}\right.
\end{equation}
We suppose that the initial condition $W(x,0)=W_0(x)$ for the
probability distribution is nonrandom. This implies that the
initial condition for $W_1(x,t)$ at $t=0$ is zero: $W_1(x,0)=
\langle \alpha(0)\widetilde W(x,0)\rangle=\langle \alpha(0)\rangle
W_0(x)=0$. The probability distribution $W(x,t)$ should be
nonnegative and normalized for all $t$: $W(x,t)\geq 0$,
$\int_{-\infty}^{\infty}W(x,t)\, dx \equiv 1$.


{\it Example 2}. Let us consider the following dynamical system
driven by two statistically independent Markovian dichotomous
noises $\aa(t)$ and $\beta(t)$:
\begin{equation}\label{Ex2}
\dot{x} = p(x)+\alpha(t) q(x)+\beta(t) g(x).
\end{equation}
The averaged probability density $W(x,t)=\langle \widetilde
W(x,t)\rangle|_{\aa,\beta}$ satisfies the following system of
master equations:
\begin{equation}\label{WEx2}
   \left\lbrace  \begin{array}{l}
 W_t + \left(p(x)W\right)_x  + \left(q(x)W_1\right)_x  + \left(g(x)P\right)_x
  =0, \\[0.5em]
 (W_1)_t + 2\nu W_1 +\left(p(x)W_1\right)_x  + \aa^2\left(q(x)W\right)_x
 +\left(g(x)Q\right)_x=0,\\[0.5em]
 (P)_t + 2\mu P +\left(p(x)P\right)_x  + \left(q(x)Q\right)_x
 +\beta^2\left(g(x)W\right)_x=0,\\[0.5em]
 (Q)_t + 2(\mu+\nu) Q +\left(p(x)Q\right)_x  + \aa^2\left(q(x)P\right)_x
 +\beta^2\left(g(x)W_1\right)_x=0,\\[0.5em]
 \end{array}\right.
\end{equation}
where the auxiliary functions $W_1(x,t)$, $P(x,t)$ and $Q(x,t)$
are some averages over realizations of the noises $\aa(t)$,
$\beta(t)$. They play the same auxiliary role as the function
$W_1(x,t)$ in the system (\ref{WW}). We suppose that $\aa(t)=\pm
\aa$, $\beta(t)=\pm \beta$ for any $t$. The characteristic
switching frequencies of these random noises are $2\nu$ and $2\mu$
respectively. We again suppose that the initial condition
$W(x,0)=W_0(x)$ is nonrandom, so the initial conditions for
$W_1(x,t)$, $P(x,t)$ and $Q(x,t)$ at $t=0$ are zeros. The
probability distribution $W(x,t)$ should be nonnegative and
normalized for all $t$: $W(x,t)\geq 0$,
$\int_{-\infty}^{\infty}W(x,t)\, dx \equiv 1$.

 {\it Example 3}. We can also consider nonlinear dynamical systems
of higher order:
\begin{equation}\label{Ex3}
   \left\lbrace  \begin{array}{l}
 \dot{x} = f(x,y)+\alpha(t) q(x,y), \\
  \dot{y} = d(x,y)+\beta(t) s(x,y),
 \end{array}\right.
\end{equation}
where $f$, $g$, $d$, $s$ are given functions. We use the same
conventions for  the noises $\aa(t)$, $\beta(t)$ and the averaged
probability distributions as in Example~2. The master equations
for the main average $P(x,y,t)$ and auxiliary averages
$P_1(x,y,t)$, $Q(x,y,t)$ and $Q_1(x,y,t)$  become
\begin{equation}\label{WEx3}
   \left\lbrace  \begin{array}{l}
 P_t + \left(fP\right)_x  + \left(gP_1\right)_x  + \left(dP\right)_y
   + \left(sQ\right)_y=0, \\[0.5em]
(P_1)_t + 2\nu P_1+ \left(fP_1\right)_x  + \aa^2\left(gP\right)_x
  + \left(dP_1\right)_y  + \left(sQ_1\right)_y=0, \\[0.5em]
 Q_t +  2\mu Q + \left(fQ\right)_x  + \left(gQ_1\right)_x  + \left(dQ\right)_y
   + \beta^2\left(sP\right)_y=0, \\[0.5em]
(Q_1)_t + 2(\mu+\nu) Q_1+ \left(fQ_1\right)_x  +
 \aa^2\left(gQ\right)_x
  + \left(dQ_1\right)_y  + \beta^2\left(sP_1\right)_y=0,
 \end{array}\right.
\end{equation}
The Cauchy initial value problem is formulated in the same way as
above: $P(x,y,0)=P_0(x,y)$ is nonrandom; $P_1(x,y,t)$, $Q(x,y,t)$
and $Q_1(x,y,t)$ at $t=0$ are zeros; $P(x,y,t)$ should be
nonnegative and normalized for all $t$: $P(x,y,t)\geq 0$,\\
$\int_{-\infty}^{\infty}\int_{-\infty}^{\infty}P(x,y,t)\, dx\, dy
\equiv 1$.

Numerous publications (see \cite{HL84,Kampen} for the
bibliography) are devoted to solution of the system (\ref{WW})
asymptotically for $t \rightarrow \infty$, that is to stationary
solutions. A number of important phenomena of noise-induced phase
transitions with applications in physics, chemistry and biology
were discovered in this stationary case.

There are only a few publications  dealing with non-stationary
solutions of (\ref{WW}). We refer to \cite{S84} and the recent
paper \cite{BB01}, where some (incomplete) exact solutions of the
system (\ref{WW}) for particular forms of the functions $p(x)$ and
$q(x)$ were obtained.

As we show in this paper, some interesting non-stationary kinetic
equations (master equations) for probability distributions allow
{\em complete} explicit closed-form solution of the general Cauchy
initial value problem. These complete solutions are obtainable
through a modification of the classical {\em Laplace cascade
method} (see e.g.\ \cite{dar-lec,forsyth,gour-l}). This method is
applicable to hyperbolic systems with two first-order linear PDEs
in the plane (as (\ref{WW}) above) or a single second-order linear
PDE in the plane. A preliminary closed-form complete solution  for
(\ref{WW}) was obtained by this method in \cite{ts05}. In that
paper a much more general method of explicit integration,
applicable to arbitrary hyperbolic higher-order linear systems (or
a single higher-order linear PDE) in the plane was developed.
Later another generalization was proposed in \cite{ts06}, it gives
closed-form complete solutions for some special class of
second-order linear hyperbolic equations with more than two
independent variables.

We give a brief account of the classical Laplace method as well as
its new generalizations in Section~\ref{sec2}. Section~\ref{sec3}
is devoted to a detailed study of the system (\ref{WW}) for the
simplest case of polynomial coefficients $p(x)=p_1x+p_2x^2$,
$q(x)=q_2x^2$, $p_1>0$, $p_2<0$,  $|p_2|>q_2>0$ (Verhulst model).
We show that for an infinite sequence of values of the switching
frequency $\nu\equiv p_1$, $\nu\equiv 2p_1$, $\nu\equiv 3p_1$,
\ldots, the complete explicit solution of the Cauchy problem is
obtainable by our methods.

In the final Section~\ref{sec4} we discuss  future prospects
and possible applications of our methods to more complicated
systems of type (\ref{WEx2}), (\ref{WEx3}).

\section{Explicit integration of hyperbolic systems}\label{sec2}

\subsection{Laplace cascade method}\label{sec21}
We give here only a special form of this method suitable for our
purpose, see \cite{dar-lec,forsyth,gour-l,ts05} for more details.

Suppose we are given a $2\times2$ first-order linear system of
PDEs
\begin{equation}\label{7}
{\left(\!\! \begin{array}{l}
 v_1 \\ v_2
\end{array}\!\!\right)\!\!}_x  \!\!=
\left( \!\!\begin{array}{ll}
 a_{11} & a_{12} \\
  a_{21} & a_{22}\!\!\!\!
\end{array}\right)
{\left( \!\!\begin{array}{l}
 v_1 \\ v_2
\end{array}\!\!\right)\!\!}_y \!\!+
\left( \!\!\begin{array}{ll}
 b_{11} & b_{12} \\
  b_{21} & b_{22}
\end{array}\!\!\right)
\left( \!\!\begin{array}{l}
 v_1 \\ v_2
\end{array}\!\!\right)
\end{equation}
with $a_{ij}=a_{ij}(x,y)$, $b_{ij}=b_{ij}(x,y)$. We will suppose
hereafter that (\ref{7}) is strictly hyperbolic, i.e.\ the
eigenvalues $\lambda_{1}(x,y)$, $\lambda_{2}(x,y)$ of the matrix
$(a_{ij})$ are real and distinct. Let $\vec p_1 = (p_{11}(x,y),
p_{12}(x,y))$, $\vec p_2 = (p_{21}(x,y), p_{22}(x,y))$ be the
corresponding left eigenvectors: $\sum_k p_{ik}a_{kj} =
\lambda_ip_{ij}$. Form the following first-order differential
operators $\X_i = \D_x - \lambda_i\D_y$ (the characteristic vector
fields) and change the initial unknown functions $v_i$ to new
characteristic functions $u_i=\sum_kp_{ik}v_k$. Then
$\X_iu_i=\sum_k (\X_i p_{ik}) v_k
 + \sum_kp_{ik}\left( (v_k)_x - \lambda_i (v_k)_y\right)=
  \sum_{k,s}p_{ik} \left(a_{ks} - \lambda_i \delta_{ks}\right)(v_s)_y +
$ \\ $ \sum_{k,s}p_{ik}b_{ks}v_s \! + \sum_k (\X_i p_{ik}) v_k
\!\!= \!\!
 \sum_s v_s \left(\sum_k p_{ik}b_{ks}\!+ (\X_i p_{is})\right)
 = \sum_k u_k \alpha_{ik}(x,y)$,
so we obtain the following {\em characteristic form} of the system
(\ref{7}):
\begin{equation}\label{6}
{\cal H}: \quad\left\lbrace  \begin{array}{l}
 \X_1u_1=\alpha_{11}(x,y)\,u_1 +\alpha_{12}(x,y)\,u_2, \\
  \X_2u_2 =\alpha_{21}(x,y)\,u_1 +\alpha_{22}(x,y)\,u_2,
\end{array}\right.
\end{equation}
The characteristic system (\ref{6}), equivalent to (\ref{7}), is
determined uniquely up to {\em operator rescaling} $\X_i
\rightarrow \gamma_i(x,y)\X_i$ and {\em gauge transformations}
$u_i \rightarrow g_i(x,y)u_i$. It is easy to check that the gauge
transformations to not change the {\em Laplace invariants} of the
system
 $h=
 \X_2(\alpha_{11}) - \X_1(\alpha_{22})
 -\X_1\X_2\ln(\alpha_{12})  -\X_1(P) +P\alpha_{11} +\alpha_{12}\alpha_{21}
 +(\alpha_{22}+\X_2(\ln\alpha_{12})+P)Q
 $ and
$k=\alpha_{12}\alpha_{21}$.
 Here $P(x,y)$ and $Q(x,y)$ are the coefficients of the commutator
\begin{equation}\label{cl}
  [\X_1,\X_2] = \X_1\X_2- \X_2\X_1 = P(x,y)\X_1 +Q(x,y)\X_2.
\end{equation}

 These invariants $h(x,y)$ and $k(x,y)$ are just the classical Laplace invariants
(cf. \cite{dar-lec,forsyth,gour-l,ts05}) of the second-order
scalar equation, obtained after elimination of $u_2$ from
(\ref{6}). Rescaling transformations of $\X_i$ change the Laplace
invariants multiplicatively: $h \rightarrow \gamma_1\gamma_2 h$,
$k \rightarrow \gamma_1\gamma_2 k$.

Starting from (\ref{6}) one can obtain {\em two} different
(inequivalent w.r.t.\ gauge transformations) second-order scalar
equations, eliminating either $u_1$ or $u_2$ from (\ref{6}).
 This observation gives rise to the
{\em Laplace cascade method of integration} of strictly hyperbolic
systems in characteristic form (\ref{6}):

$({\cal L}_1)$ If $k$ vanishes then either $\aa_{12}$ or
$\aa_{21}$ vanishes so the system becomes triangular:
\begin{equation}\label{Tr}
\left\lbrace  \begin{array}{l}
 \X_1u_1=\alpha_{11}\,u_1,  \\
  \X_2u_2 =\alpha_{21}\,u_1 +\alpha_{22}\,u_2,
\end{array}\right.
{\textrm{or}} \quad \left\lbrace  \begin{array}{l}
 \X_1u_1=\alpha_{11}\,u_1 +\alpha_{12}\,u_2, \\
  \X_2u_2 =\alpha_{22}\,u_2.
\end{array}\right.
\end{equation}
If we perform an appropriate change of coordinates $(x,y)
\rightarrow (\overline x, \overline y)$ (NOTE: for this we have to
solve  first-order nonlinear ODEs $dy/dx = -\lambda_i(x,y)$, cf.\
Appendix in \cite{grig2}) one can suppose $\X_1=\D_{\overline x}$,
$\X_2=\D_{\overline y}$ and obtain the complete solution of
(\ref{Tr}) in quadratures: if for example $\aa_{12} \equiv 0$,
then
\begin{equation}\label{solu12}
  \begin{array}{l}
 u_1=Y(\overline y)\exp\left(-\int \alpha_{11} \, d\overline x \right),
 \\[1ex]
 u_2 =\exp\left(-\int \alpha_{22} \, d\overline y \right)
  \left( X(\overline x) + \int Y(\overline y)
\exp\left(\int (\alpha_{11}\, d\overline x - \alpha_{22}\,
d\overline y)\right) d\overline y \right)
\end{array} 
\end{equation}
where $ X(\overline x)$ and $Y(\overline y)$ are two arbitrary
functions of the characteristic variables $\overline x$,
$\overline y$ respectively.

$({\cal L}_2)$ If $k \neq 0$, transform the system into  a
second-order scalar equation eliminating $u_2$ from (\ref{6}):
from the first equation
\begin{equation}\label{difsub}
u_2 = (\X_1 u_1 -\alpha_{11} u_1)/\aa_{12},
\end{equation}
substitute this expressions into the second equation obtaining
$\hatl u_1=\X_2\frac{1}{\aa_{12}}(\X_1u_1 - \aa_{11}u_1) -
\aa_{21}u_1 - \frac{\aa_{22}}{\aa_{12}}(\X_1u_1 - \aa_{11}u_1) =
0$. Now, using the commutator relation (\ref{cl}), we can
represent $\hatl u_1$ as $\hatl u_1 = (\X_1\X_2 + \beta_1\X_1 +
\beta_2\X_2 +  \beta_3)u_1 = (\X_1+\beta_2) (\X_2+\beta_1)u_1-hu_1
= 0$. From this form we see that this equation is equivalent to
another $2 \times 2$ system
\begin{equation}\label{H1}
{\cal H}_{(1)}\ : \quad\left\lbrace
  \begin{array}{l}
  \X_2u_1=-\beta_1u_1 +\overline u_2, \\
  \X_1\overline u_2 =hu_1-\beta_2\overline u_2.
\end{array}\right.
\end{equation}
This new system (we will call it {\em $X_1$-transformed} system)
has the same characteristic form (\ref{6}) with {\em different}
coefficients in the right-hand side. It also has new Laplace
invariants $ h_{(1)}$, $k_{(1)}$, and it turns out that $k_{(1)}$
equals to the invariant $h$ of the original system. So if we have
$k_{(1)}=h=0$,
 we solve this new system in quadratures and
using the same differential substitution (\ref{difsub}) we obtain
the complete solution of the original equation $\hatl u=0$.

$({\cal L}_3)$ If again $k_{(1)} \neq 0$, apply this
$X_1$-transformation several times, obtaining a sequence of
$2\times 2$ characteristic systems ${\cal H}_{(2)}$, ${\cal
H}_{(3)}$, \ldots\ \
 If on any step we get $k_{(m)}=0$,
we solve the corresponding system in quadratures and, using the
differential substitutions (\ref{difsub}), obtain the complete
solution of the original system. Alternatively one may perform
{\em $\X_2$-transformations}, eliminating $u_1$ instead of $u_2$
on step $({\cal L}_2)$. In fact this $\X_2$-transformation is a
reverse of the $\X_1$-trans\-for\-ma\-tion up to a gauge
transformation (see \cite{anderson}). So we have (infinite in
general) chain of systems
\begin{equation}\label{ch}
   \ldots \stackrel{\X_2}{\leftarrow} {\cal H}_{(-2)}
    \stackrel{\X_2}{\leftarrow}
   {\cal H}_{(-1)}\stackrel{\X_2}{\leftarrow}   {\cal H}
    \stackrel{\X_1}{\rightarrow}
    {\cal H}_{(1)} \stackrel{\X_1}{\rightarrow} {\cal H}_{(2)}
    \stackrel{\X_1}{\rightarrow}
     \ldots
\end{equation}
and the corresponding chain of their Laplace invariants
\begin{equation}\label{li}
   \ldots , k_{(-3)},\  k_{(-2)},\  k_{(-1)},\  k,\
    k_{(1)} = h,\ k_{(2)}, \ k_{(3)}, \ldots
\end{equation}
 We do not need to
keep the invariants $h_{(i)}$ in (\ref{li}) since
$k_{(i)}=h_{(i-1)}$. If on any step we have $k_{(N)}=0$ then the
chains (\ref{ch}) and (\ref{li}) can not be continued: the
differential substitution (\ref{difsub}) is not defined; precisely
on this step the corresponding system (\ref{6}) is triangular and
we can find its complete solution as well as the complete solution
for any of the systems of the chain (\ref{ch}).

 As one may prove (see e.g.\ \cite{dar-lec}) if the chain
(\ref{ch}) is finite in both directions (i.e.\ we have
$k_{(N)}=0$, $k_{(-K)}=0$ for some $N\geq 0$, $K\geq 0$) one may
even obtain a quadrature-free expression for the general solution
of the original system:
\begin{equation}\label{XY}
  \begin{array}{l}
    u_1 = c_0F + c_1F' + \!\ldots \!+ \! c_NF^{(N)}\!\! +
d_0\widetilde G + d_1\widetilde G' +
   \ldots +  d_{K+1}\widetilde G^{(K+1)},\\[1ex]
    u_2 = e_0F + e_1F' + \!\ldots \!+ \! e_NF^{(N)}\!\! +
f_0\widetilde G + f_1\widetilde G' +
   \ldots +  f_{K+1}\widetilde G^{(K+1)},
\end{array}
\end{equation}
with definite $c_i(\overline x,\overline y)$, $d_i(\overline
x,\overline y)$, $e_i(\overline x,\overline y)$, $f_i(\overline
x,\overline y)$ and
 two arbitrary functions  $F(\overline x)$, $\widetilde G(\overline y)$
of the characteristic variables. Vice versa: existence of ({\em a
priori} not complete) solution of the form (\ref{XY}) with
arbitrary functions $F$, $G$ of characteristic variables implies
$k_{(s)}=0$, $k_{(-r)}=0$ for some $s \leq N$, $r \leq K$. So {\em
minimal differential complexity} of the answer (\ref{XY}) (the
number of terms in it) is equal to the number of steps necessary
to obtain vanishing Laplace invariants in the chains (\ref{ch}),
(\ref{li}) and consequently triangular systems. Complete proofs of
these statement may be found in \cite[t.~2]{dar-lec},
\cite{forsyth,gour-l} for the case $\X_1=\D_{x}$, $\X_2=\D_{y}$,
for the general case cf.\ \cite[p. 30]{gour-l} and
\cite{anderson}.

We give a detailed example of application of this method in
Section~\ref{sec3}.

There were some attempts to generalize Laplace transformations for
higher-order systems or the number of independent variables larger
than 2, both in the classical time \cite{LeRoux,pisati,Petren} and
in the last decade \cite{Athorne,ts05}. As one can show, all of
them essentially try to triangulize the given system in some
sense. A general definition of ``generalized factorization''
(triangulation) comprising all known practical methods was given
in \cite{ts98}. Unfortunately the theoretical considerations of
\cite{ts98} did not provide any algorithmic way of establishing
generalized factorizability of a given higher-order operator or a
given higher-order system. Below we present other approach for
search of ``generalized factorizations'' resulting in explicit
complete solution of some classes of hyperbolic systems.

\subsection{Generalized Laplace cascade method for $n \times n$
hyperbolic systems in the plane }\label{sec22}

Any $n\times n$ first-order linear system
\begin{equation}\label{3.3}
(v_i)_x = \sum_{k=1}^n a_{ik}(x,y)(v_k)_y +
  \sum_{k=1}^n b_{ik}(x,y)v_k
\end{equation}
with strictly hyperbolic matrix $(a_{ik})$ (i.e.\ with real and
distinct eigenvalues of this matrix) is equivalent to a system in
characteristic form 
\begin{equation}\label{3.2}
 \X_iu_i= \sum_k \alpha_{ik}(x,y)u_k,
\end{equation}
as a straightforward calculation similar to that in the beginning
of Section~\ref{sec21} immediately shows.

Our generalization of the Laplace transformations consists in the
following.

$({\cal L}_1)$ For a given $n\times n$ characteristic system
(\ref{3.2}) choose one of its equations with a non-vanishing
off-diagonal coefficient  $\alpha_{ik}\neq 0$, find
 $u_k = (\X_iu_i - \sum_{s \neq k} \alpha_{is}u_s)/\alpha_{ik}$
and substitute this expression into all other equations of the
system. We obtain one second-order equation
\begin{equation}\label{L11}
   \X_k \big(\frac{1}{\!\alpha_{ik}\!}(\X_iu_i - \!\sum_{s \neq k}\alpha_{is}u_s)\big)
    - \!\sum_{p\neq k}\alpha_{kp}u_p
  - \frac{\alpha_{kk}}{\alpha_{ik}}(\X_iu_i - \!\sum_{s \neq k}\alpha_{is}u_s)\! =\!
  0
\end{equation}
and $n-2$ first-order equations
\begin{equation}\label{L12}
  \X_j u_j -\sum_{s \neq k}\alpha_{js}u_s -
  \frac{\alpha_{jk}}{\alpha_{ik}}(\X_iu_i - \sum_{s \neq k}\alpha_{is}u_s) =
  0
\end{equation}
for $j \neq i,k$.

$({\cal L}_2)$  The second step consists in rewriting the system
(\ref{L11}), (\ref{L12}) in the following form with slightly
modified unknown functions $\overline u_j = u_j + \rho_j(x,y)u_i$,
$j \neq i,k$, $\overline u_i \equiv u_i$, new coefficients
$\beta_{pq}(x,y)$ but the same characteristic operators $\X_p$:
\begin{equation}\label{L21}
\!\!  \X_i (\frac{1}{\alpha_{ik}\!\!}(\X_k\overline u_i
   - \!\sum_{s \neq k}\beta_{is}\overline u_s))
    - \!\sum_{p\neq k}\beta_{kp}\overline u_p
  - \frac{\beta_{kk}\!}{\alpha_{ik}\!}(\X_k\overline u_i -
    \!\sum_{s \neq k}\beta_{is}\overline u_s) \!=\!
  0,
\end{equation}
\begin{equation}\label{L22}
 \X_j \overline u_j -\sum_{s \neq k}\beta_{js}\overline u_s -
  \frac{\beta_{jk}}{\alpha_{ik}}(\X_k\overline u_i - \sum_{s \neq k}\beta_{is}\overline u_s) =
  0, \quad j \neq i,k.
\end{equation}
As one can prove this is always possible in a unique way.


$({\cal L}_3)$  Introducing $\overline u_k
=\frac{1}{\alpha_{ik}}(\X_k\overline u_i - \sum_{s \neq
k}\beta_{is}\overline u_s)$ rewrite (\ref{L21}), (\ref{L22}) as
the transformed characteristic system
\begin{equation}\label{L31}
      \left\lbrace  \begin{array}{l}
 \X_i \overline u_k =  \sum_{p}\beta_{kp}\overline u_p ,  \\[0.5em]
 \X_k\overline u_i = \sum_{s \neq k}\beta_{is}\overline u_s
  +\alpha_{ik}\overline u_k, \\[0.5em]
 \X_j \overline u_j =\sum_{s}\beta_{js}\overline u_s,
  \quad j \neq i,k.
 \end{array}\right.
\end{equation}

The reason of doing such generalized Laplace transformation
consists in the fact that after it (or after a chain of such
transformations) one may obtain a triangular system (\ref{3.2}),
solve it in quadratures and doing the inverse steps $({\cal
L}_3)$, $({\cal L}_2)$, $({\cal L}_1)$ obtain the complete
solution of the original $n \times n$ hyperbolic system.

Cf.~ \cite{ts05} for the details, an example and the proof of
correctness of the step $({\cal L}_2)$.

\subsection{Explicit solution of equations with more than two
independent variables}\label{sec23}

The method described in this Section is based on an idea given by
Ulisse Dini in 1902 (cf.~\cite{ts06} for references and details).
We will limit here to a simple example showing the idea of the
method.

Let us take the following equation:
\begin{equation}\label{dex}
 Lu = (\hat D_x\hat D_y + x \hat D_x\hat D_z - \hat D_z)u =0.
\end{equation}
It has three independent derivatives $\hat D_x$, $\hat D_y$, $\hat
D_z$, so the Laplace method is \emph{not} applicable. On the other
hand its principal symbol splits into product of two first-order
factors: $\xi_1\xi_2 + x \xi_1\xi_3 =\xi_1(\xi_2+x\xi_3)$. This is
no longer a typical case for hyperbolic operators in
dimension~$3$; we will use this special feature introducing two
characteristic operators $\hat X_1=\hat D_x$, $\hat X_2=\hat D_y +
x \hat D_z$. We have again a nontrivial commutator  $[\hat
X_1,\hat X_2] =  \hat D_z= \hat X_3$. The three operators $\hat
X_i$ span the complete tangent space in every point $(x,y,z)$.
Using them one can represent the original second-order operator in
one of two partially factorized forms:
$$ L = \hat X_2\hat X_1 - \hat X_3 =  \hat X_1\hat X_2 - 2\hat X_3.$$
Let us use the first one and transform the equation into a system
of two first-order equations:
\begin{equation}\label{Dini2e}
 Lu=0 \Longleftrightarrow
   \left\lbrace  \begin{array}{l}
\hat X_1 u = v, \\ \hat X_3 u = \hat X_2 v.
\end{array}\right.
\end{equation}
Here comes the difference with the classical case $dim=2$: we can
not express $u$ as we did in (\ref{difsub}). But we have another
obvious possibility instead: cross-differentiating the left hand
sides of (\ref{Dini2e}) and using the obvious identity $[\hat
X_1,\hat X_3] = [ \hat D_x, \hat D_z]=0$ we get $  \hat X_1 \hat
X_2v =   \hat D_x (\hat D_y + x\hat D_z) v = \hat X_3 v=\hat D_z v
$ or $ 0=\hat D_x (\hat D_y + x\hat D_z) v - \hat D_z v = (\hat
D_x \hat D_y + x \hat D_x \hat D_z) v
 = (\hat D_y + x\hat D_z)  \hat D_x v = \hat X_2\hat X_1 v$.

Since we have now  another second-order equation which is
``naively'' factorizable we easily find its complete solution:
$$v= \int \phi(x,xy-z) \, dx + \psi(y,z)$$
where $\phi$ and $\psi$ are two arbitrary functions of two
variables each; they give the general solutions of the equations
$\hat X_2\phi=0$, $\hat X_1\psi=0$.

Now we can find  $u$:
$$ u= \int \Big(v\, dx  + (\hat D_y + x\hat D_z)v\, dz \Big)+ \theta(y),
$$
where an extra free function $\theta$ of one variable appears as a
result of integration in (\ref{Dini2e}).

So we have seen that such {\em Dini transformations}
(\ref{Dini2e}) in some cases may produce a complete solution in
explicit form for a non-trivial three-dimensional equation
(\ref{dex}). This explicit solution can be used to solve initial
value problems for (\ref{dex}).

\section{Verhulst model}\label{sec3}

Here we describe in detail the procedure of solution for the
system~(\ref{WW}).

The characteristic operators and left eigenvectors of this $2
\times 2$-system are simple: $\X_i=\D_t-\lambda_i\D_x$,
$\lambda_{1,2}= -p(x) \pm q(x)$, $p_{11}=p_{21}=p_{22}=1$,
$p_{12}=-1$. The characteristic system (\ref{6}) for the new
characteristic functions $u_1=W-W_1$, $u_2=W+W_1$ is
\begin{equation}\label{WWch}
   \left\lbrace  \begin{array}{l}
 \X_1u_1= -(p_x-q_x+\nu)\,u_1 +\nu \,u_2  ,\\
 \X_2u_2 = \nu \,u_1 - (p_x+q_x+\nu)\,u_2.
\end{array}\right.
\end{equation}
The Laplace invariants are
 $h= \nu^2 - [ p_{xx}q^2(p+q) + p_x^2q^2 - p_xq_xq(3p+q)
      -q_{xx}pq(p+q) - q_x^2p(2p+q) ]/q^2 $,
 $k=\nu^2$, so if $\nu$,
$p(x)$ and $q(x)$ satisfy a second-order differential relation
$h=0$, one can solve (\ref{WW}) in quadratures. Especially simple
formulas may be obtained for polynomial $p(x)=p_1x+p_2x^2$,
$q(x)=q_2x^2$: in this case $k=\nu^2$,
$h=h_{(-2)}=\nu^2-p_1^2$ so if $\nu=p_1$, one may solve (\ref{WW})
explicitly. It is convenient at this point to use the
dimensionless variable $\tau=\nu t$; so we have to change $t
\mapsto \tau$, $\nu\mapsto 1$, $p_1\mapsto 1$ and change $p_2$,
$q_2$ respectively. For simplicity we will still use the same
notations $p_2$, $q_2$.

After the necessary transformation, described in
Section~\ref{sec22}, we obtain the following quadrature-free
expression for the complete solution of the system (\ref{WW}):
\begin{equation}\label{WFG}
\begin{array}{l}
\displaystyle W= \frac{q_2}{x^2}\left[ F'(\overline x)
-F(\overline x) +  G'(\overline y) -  G(\overline y)
\right],\\[0.5em]
\begin{array}{l} \displaystyle  \!\! W_1=
\frac{1}{x^3}\left[ - q_2xG'(\overline y) +(1+p_2x)G(\overline y)
 +  q_2 xF'(\overline x)\right. \\[0.5em]
 \quad \left. {} + (1+p_2x) F(\overline x) \right],
\end{array}
\end{array}
\end{equation}
where $\overline x = -t + \ln\frac{x}{1+(p_2+q_2)x}$, $\overline y
= -t + \ln\frac{x}{1+(p_2-q_2)x}$ are the characteristic variables
($\X_2\overline x =0$, $\X_1\overline y =0$) and $F$, $G$ are two
arbitrary functions of the corresponding characteristic variables.

For the case $\nu^2 \!\neq \!p_1^2$ we can compute other Laplace
invariants of the chain (\ref{li}):
 $h_{(1)}\!=h_{(-3)}\!=\nu^2-4p_1^2$,
 $h_{(2)}=h_{(-4)}=\nu^2-9p_1^2$,
 $h_{(3)}=h_{(-5)}=\nu^2-16p_1^2$, etc.,
so for the fixed $p(x)=p_1x+p_2x^2$, $q(x)=q_2x^2$ and $\nu=\pm
p_1$, $\nu=\pm 2p_1$, $\nu=\pm 3p_1$, \ldots\ one can obtain
closed-form quadrature-free complete solution of the system
(\ref{WW}), with increasing complexity of the answer (\ref{XY}).

Now we demonstrate how the formulas (\ref{WFG}) may be used to
solve the Cauchy initial value problem. For this set $\tau=0$
inside the variables $\overline x$, $\overline y$ and equate
$W(x,0)=W_0(x)$, $W_1(x,0)=0$. Since now $\overline x =
\ln\frac{x}{1+(p_2+q_2)x}$, $\overline y =
\ln\frac{x}{1+(p_2-q_2)x}$, one can express the derivatives $F' =
dF/d\overline x$,  $G' = dG/d\overline y$ as $F' =
\frac{dF}{dx}\frac{dx}{d\overline x}$, $G' =
\frac{dG}{dx}\frac{dx}{d\overline y}$ and obtain from (\ref{WFG})
a system of two linear ODEs for $F(x)=F(\overline x(x))$,
$G(x)=G(\overline y(x))$. It may be solved explicitly (see an
explanation of this fact in Section~\ref{sec4}) for any $W_0(x)$:
\begin{equation}\label{nFG}
\begin{array}{l}
  F(x)=\frac{1}{2q_2^2(1+(p_2+q_2)x)}
   \left[-x \int_{c_0}^x\frac{W_0(\theta)}{\theta} d\theta
     + (1+q_2x)\int_{c_1}^x W_0(\theta)\, d\theta \right],\\[2ex]
  G(x)=\frac{1}{2q_2^2(1+(p_2-q_2)x)}
   \left[x \int_{c_0}^x\frac{W_0(\theta)}{\theta} d\theta
     + (q_2x-1)\int_{c_1}^x W_0(\theta)\, d\theta \right].
\end{array}
\end{equation}
Perform now the inverse substitution $F(\overline x)=F(x(\overline
x))$, $G(\overline y)=G(x(\overline y))$ (for $\tau=0$) to find
the ``true'' functions $F(\overline x)$, $G(\overline y)$ suitable
for substitution into (\ref{WFG}) for any $\tau$. This final form
of the explicit solution of the Cauchy problem is:
$$
W(x,\tau)= \frac{1}{2q_2x^2}\left[I_1(\hat y) - I_1(\hat x)
\right]
 + \frac{W_0(\hat y)}{2(e^\tau(1+(p_2-q_2)x) - x(p_2-q_2))^2}
$$
$$
{} + \frac{W_0(\hat x)}{2(e^\tau(1+(p_2+q_2)x) - x(p_2+q_2))^2},
$$

\bigskip

$$
W_1(x,\tau)=\frac{I_1(\hat y)- I_1(\hat
x)}{2q_2^2x^3}\left[(e^{-\tau}-1)p_2x -1 \right]+
\frac{e^{-\tau}}{2q_2^2x^2}\left[I_2(\hat y) - I_2(\hat x) \right]
$$
$$
 - \frac{W_0(\hat y)}{2(e^\tau(1+(p_2-q_2)x) - x(p_2-q_2))^2}
 + \frac{W_0(\hat x)}{2(e^\tau(1+(p_2+q_2)x) - x(p_2+q_2))^2},
$$

\medskip

where $\hat x = \frac{x}{(e^\tau(1+(p_2+q_2)x) - x(p_2+q_2))}$,
$\hat y = \frac{x}{(e^\tau(1+(p_2-q_2)x) - x(p_2-q_2))}$,
 $I_1(z)= \int_{c_1}^z W_0(\theta)\, d\theta$,
 $I_2(z)= \int_{c_0}^z \frac{W_0(\theta}{\theta} d\theta$, $c_0$
and $c_1$ may be chosen arbitrary.

One can check that $\hat x < \hat y$ for all $t\geq 0$, $x\geq 0$.

We get an especially simple form of this solution for the initial
distribution
 $W_0(x)=\delta(x-x_*)$ with some fixed initial state $x(0)=x_*>0$:
$$
W(x,\tau)= \frac{\delta(\hat x-x_*)}{2(e^\tau(1+(p_2+q_2)x) -
x(p_2+q_2))^2}
 + \frac{\delta(\hat y-x_*)}{2(e^\tau(1+(p_2-q_2)x) -
x(p_2-q_2))^2}
$$

\medskip

$$
{} + \frac{H(\hat y-x_*) - H(\hat x-x_*)}{2q_2x^2}.
$$
Here $H(z)= \int_{-\infty}^z \delta(\theta)\, d\theta$ is the
Heaviside function.

According to the standard formula
$\delta(\phi(x))=\delta(\phi^{-1}(0))/\phi'(\phi^{-1}(0))$ one
gets $\displaystyle \frac{\delta(\hat
x-x_*)}{2(e^\tau(1+(p_2+q_2)x) - x(p_2+q_2))^2} = \frac{\delta(
x-\frac{e^{\tau}x_*}{1-(p_2+q_2)(e^\tau -1)x_*})}{2e^\tau}$,

$\displaystyle \frac{\delta(\hat y-x_*)}{2(e^\tau(1+(p_2-q_2)x) -
x(p_2-q_2))^2} = \frac{\delta(
x-\frac{e^{\tau}x_*}{1-(p_2-q_2)(e^\tau -1)x_*})}{2e^\tau}$, so we
see that this simple solution (and consequently the complete
solution) obviously obeys the necessary physical requirements of
positivity and normalization:  $W(x,t)\geq 0$,
$\int_{-\infty}^{\infty}W(x,t)\, dx \equiv 1$. Asymptotically, for
$\tau \rightarrow \infty$, this solution exponentially fast
converges to the stationary probability distribution
$W_\infty(x)=0$ outside the interval $\frac{1}{|p_2-q_2|}< x <
\frac{1}{|p_2+q_2|}$ and $W_\infty(x)=1/(2q_2x^2)$ inside this
interval.

\section{Concluding remarks and future prospects}\label{sec4}

There is an algorithmic possibility to obtain closed-form
solutions of the Cauchy problem for the more complicated cases
$\nu=mp_1$, $m=2,3,4,\ldots$ in the Verhulst model. The respective
classical form (\ref{XY}) is algorithmically obtainable with the
methods of Section~\ref{sec21}. Since the orders of derivations of
$F(\overline x)$, $G(\overline y)$ in the right-hand sides of
(\ref{XY}) are proportional to the integer coefficient $m$ in the
relation $\nu=mp_1$, directly assigning $W(x,0)=W_0(x)$,
$W_1(x,0)=0$ in this formula for $\tau=0$ will result in a linear
system of ODEs for $F(\overline x)$, $G(\overline y)$ of high
order with nonconstant coefficients. Much more efficient is to use
the transformations (\ref{difsub}) directly: simply recalculate
the Cauchy data for the new functions $\overline u_2$ on step
$({\cal L}_2)$, using (\ref{H1}), until we get (after $m$ steps) a
triangular system, solve this system for the recalculated Cauchy
data and then use the inverse $\X_2$-transformations to get the
solution of the original system. This also explains why we could
find the solution (\ref{nFG}) in the case $\nu=p_1$ in
Section~\ref{sec3}.

Methods, described in Sections~\ref{sec22},~\ref{sec23}, suggest
that one can also investigate systems (\ref{WEx2}), (\ref{WEx3})
and classify completely integrable cases for special forms of
their coefficients $p(x)$, $g(x)$, $ f(x,y)$, $q(x,y)$, $d(x,y)$,
$s(x,y)$ and switching frequencies $\mu$, $\nu$.Systematic
investigation of such integrable cases will be reported in
subsequent publications.

\end{document}